\documentclass[twoside]{amsart}
\usepackage{graphicx,xcolor} 
\usepackage{latexsym}
\usepackage{enumitem}
\usepackage{amssymb,amsmath,amsopn}
\usepackage{subcaption} 
\usepackage{array}
\usepackage{hyperref}
\usepackage{cite}
\usepackage{tabularx}

\renewcommand{\Re}{\mathbb R}

\newcommand{\BB}{\mathbf B}

\newcommand{\Sph}{\mathbb{S}}

\newcommand{\F}{\mathcal{F}}

\newcommand{\K}{\mathcal{K}}
\newcommand{\Q}{\mathcal{Q}}
\newcommand{\rin}{\rho^{in}}
\newcommand{\rex}{\rho^{ex}}

\DeclareMathOperator{\card}{card}
\DeclareMathOperator{\inter}{int}

\DeclareMathOperator{\conv}{conv}

\DeclareMathOperator{\perim}{perim}
\DeclareMathOperator{\area}{area}

\DeclareMathOperator{\surf}{surf}

\DeclareMathOperator{\vol}{vol}
\DeclareMathOperator{\bd}{bd}

\theoremstyle{plain}
\newtheorem{theorem}{Theorem}[section]

\theoremstyle{definition}
\newtheorem{definition}[theorem]{Definition}
\numberwithin{equation}{section}

\begin{document}

\title[Centroids and equilibrium points]{Centroids and equilibrium points of convex bodies}

\author[Z. L\'angi]{Zsolt L\'angi}
\author[P. L. V\'arkonyi]{P\'eter L. V\'arkonyi}

\address{Zsolt L\'angi, Department of Algebra and Geometry, Budapest University of Technology and Economics,\\
M\H uegyetem rkp. 3., H-1111 Budapest, Hungary, and\\
Alfr\'ed R\'enyi Institute of Mathematics,\\
Re\'altanoda utca 13-15, H-1053, Budapest, Hungary}
\email{zlangi@math.bme.hu}
\address{P\'eter L. V\'arkonyi, Department of Mechanics, Materials, and Structures, Budapest University of Technology and Economics,\\
M\H uegyetem rkp. 3., H-1111 Budapest, Hungary}
\email{varkonyi.peter@epk.bme.hu}

\thanks{The authors thank Timea Szab\'o and G\'abor Domokos for their permissions to use their images of bodies in various equilibrium classes (Fig. 1). Zsolt L\'angi is partially supported by the ERC Advanced Grant ``ERMiD'', and the NKFIH grant K147544. P. L. V\'arkonyi is supported by the NKFIH grant 143175.}


\begin{abstract}
The aim of this note is to survey the results in some geometric problems related to the centroids and the static equilibrium points of convex bodies. In particular, we collect results related to Gr\"unbaum's inequality and the Busemann-Petty centroid inequality, describe classifications of convex bodies based on equilibrium points, and investigate the location and structure of equilibrium points, their number with respect to a general reference point as well as the static equilibrium properties of convex polyhedra. 
\end{abstract}

\maketitle
 
\section{Introduction}\label{sec:intro}

The investigation of the centroids and static equilibrium points of solids has been in the focus of research since ancient Greece. These concepts play an important role in many different disciplines from physics to engineering to geology, and also appear in many problems in pure mathematics. The aim of this paper is to collect the results in the geometric literature related to some of these problems.

Apart from a few remarks, the objects of our paper are convex bodies in the $d$-dimensional Euclidean space $\Re^d$, defined as compact, convex sets with nonempty interior. We denote the standard inner product of this space by $\langle \cdot, \cdot \rangle$, and the Euclidean norm of a vector $p \in \Re^d$ by $|p| = \sqrt{\langle p,p \rangle}$. The \emph{closed segment} with endpoints $p,q \in \Re^d$ is denoted as $[p,q]$. An important convex body in $\Re^d$ is the closed unit ball centered at the origin $o$, which we denote by $\BB^d$. The boundary of this ball is denoted by $\Sph^{d-1}$. We use the standard notation $\inter (X)$, $\bd (X)$, $\conv (X)$, $\vol_d (X)$ for the \emph{interior, boundary, convex hull} and $d$-dimensional \emph{volume} of a set $X \in \Re^d$. Two-dimensional volume is also called \emph{area}, denoted by $\area(\cdot)$, and $3$-dimensional volume is denoted by $\vol(\cdot)$. For any $s > 0$, we set $\kappa_s = \frac{\pi^{s/2}}{\Gamma \left( \frac{s}{2} + 1 \right)}$, where $\Gamma(x) = \int_0^{\infty} t^{x-1} e^{-t} \, dt$ is the usual Gamma function. Then, for any nonnegative integer $d$, $\kappa_d = \vol_d (\BB^d)$.

As usual, the orthogonal complement of a linear subspace $L$ of $\Re^d$ is denoted by $L^{\perp}$, and the orthogonal projection of a set $X$ onto $L$ is denoted by $X | L$. For brevity, if $L$ is a line generated by a nonzero vector $u$, we set $u^{\perp}=L^{\perp}$, $K | u = K | L$, and $K | u^{\perp} = K | L^{\perp}$.
Finally, for any $u \in \Sph^{d-1}$ and $\alpha \in \Re$, we denote the hyperplane $\{ x \in \Re^d: \langle x, u \rangle = \alpha \}$ by $H(u,\alpha)$. In this notation $u$ is a \emph{normal vector} of $H(u,\alpha)$ and $|\alpha|$ is its distance from $o$.

Let $K$ be a convex body with $o \in \inter(K)$. The \emph{radial function} $\rho_K : \Sph^{d-1} \to \Re$ of $K$ is defined as
    \begin{equation} \label{eq:radialfunction}
        \rho_K(u)=\max\{\alpha\in\Re:\alpha u\in K\},
    \end{equation}
and its \emph{support function} $h_K : \Sph^{d-1} \to \Re$ as
    \begin{equation} \label{eq:supportfunction}
        h_K(u)=\max \{ \langle x,u\rangle : x \in K \}.
    \end{equation}
The latter function is often defined for any $u \in \Re^d$ using the formula in (\ref{eq:supportfunction}). To distinguish these two functions, we use the notation $\bar{h}_K : \Re^d \to \Re$ for the one extended to $\Re^d$. We note that according to our definition, the function $\bar{h}_K$ is a positively homogeneous, convex, nonnegative function which is zero only at $o$. It is also known \cite{Schneiderbook} that any function satisfying these properties is the support function of some convex body $K$ with $o \in \inter(K)$.

A set $S \subseteq \Re^d$ is called \emph{star-shaped} with respect to $p$, if for any $q \in S$, we have $[p,q] \subseteq S$. In particular, a convex set $K$ is star-shaped with respect to any $p \in K$. We define the radial function $\rho_S : \Sph^{d-1} \to \Re$ of any compact set, star-shaped with respect to $o$, by the formula in (\ref{eq:radialfunction}). Such a set is called a \emph{star body} if its radial function is strictly positive and continuous.


In the remaining part of the paper we present the results related to our basic concepts. In particular, in Section~\ref{sec:centroid} we investigate the centroids of convex bodies, and discuss results related to two important problems concerning them: Gr\"unbaum's centroid inequality in Subsection~\ref{subsec:Grunbaum}, and the Busemann-Petty centroid inequality in Subsection~\ref{subsec:BP}. In Section~\ref{sec:equilibria} we study various problems related to the static equilibrium points of convex bodies. Within this, in Subsection~\ref{subsec:primary} we present a classification system of convex bodies based on the number of their equilibria. In Subsection~\ref{subsec:sectert} we refine this system by distinguishing convex bodies based on the topological properties of the gradient flow of their radial functions. In Subsection~\ref{subsec:transition} we investigate how difficult it is to move a convex body from one class to another one. In Subsections~\ref{subsec:location}, \ref{subsec:inhom} and \ref{subsec:polyhedra} we deal with the locations of equilibrium points on a convex body, equilibria with respect to a general reference point, and equilibria of convex polyhedra. Finally, in Section~\ref{sec:appl}, we collect some applications and problems related to the presented topics.

\section{The centroid of a convex body}\label{sec:centroid}

The main concept of this section is the following.

\begin{definition}\label{centroid}
Let $X$ be a compact set in $\Re^d$ with nonempty interior. Then the 
\emph{centroid} of $K$ is the point
\[
c(X) = \frac{1}{\vol_d(X)} \int_{X} x \, dx.
\]
\end{definition}

In the above definition, if $c(X) = o$, we say that $X$ is \emph{centered}. A concept related to the centroid of a set is its first moment with respect to a hyperplane.

\begin{definition}\label{defn:firstmoment}
Let $X$ be a compact set in $\Re^d$ with nonempty interior, and let $u \in \Sph^{d-1}$, $\alpha \in \Re$. Then the quantity
\begin{equation}\label{eq:firstmoment}
\int_X \left( \langle u , x \rangle - \alpha \right) \, dx
\end{equation}
is called the \emph{first moment} of $K$ with respect to the hyperplane $H(u,\alpha)$.
\end{definition}

We note that the quantity $\langle u,x \rangle - \alpha$ in (\ref{eq:firstmoment}) is equal to the signed distance of $x$ from $H(u,\alpha)$, according to the orientation induced by $u$. From Definitions~\ref{centroid} and \ref{defn:firstmoment} it follows that the first moment of $X$ with respect to $H(u,\alpha)$ is equal to $\langle \vol_d(X) c(X), u \rangle - \vol_d(X) \alpha$, that is, it is equal to the first moment of a point mass at $c(X)$ with weight $\vol_d(X)$. While the notion of centroid in non-Euclidean spaces is outside the scope of this paper, we remark that an axiomatic generalization of centroid in various manifolds including hyperbolic and spherical spaces, based partly on this observation, can be found in \cite{Galperin2}. For more information on centroids in spherical or hyperbolic spaces, the interested reader is referred to \cite{Bjerre, Fog, Galperin1, BHPS2023}.

The geometric concept of a centroid is partly motivated by the mechanics of rigid bodies. In this context the centroid of a convex body $K$ coincides with the \emph{center of mass} of the solid $K$ in the case of homogeneous mass distribution, as well as with the \emph{center of gravity} of $K$ in a uniform gravitational field. However a physical object may also have inhomogeneous mass distribution, in which case the center of mass and the center of gravity of $K$ are identical to each other but different from the centroid, as they can be any interior point of $K$ depending on the mass distribution.

\subsection{Gr\"unbaum's centroid inequality}\label{subsec:Grunbaum}

Consider a convex disk $K\subset \Re^2$, and let $L$ be a line through its centroid. This line dissects $K$ into two convex disks $K_1$ and $K_2$. It was observed by Winternitz in 1923 \cite{Blaschkebook} that for any such disk $K$ and line $L$, we have
\[
\frac{4}{5} \leq \frac{\area(K_1)}{\area(K_2)} \leq \frac{5}{4},
\]
with equality if $K$ is a triangle.

This result was generalized by Gr\"unbaum \cite{Grunbaum1960} for any dimension in 1960 as follows.

\begin{theorem}\label{thm:Grunbaum}
Let $K \subset \Re^d$ be a centered convex body. For any $u \in \Sph^{d-1}$ let $K_u$ denote the intersection of $K$ with the closed half space bounded by $u^{\perp}$ and containing $u$. Then, for any $u \in \Sph^{d-1}$, we have
\[
\frac{\vol_d(K_u)}{\vol_d(K_{-u})} \leq \left( 1 + \frac{1}{d} \right)^{d} -1.
\]
\end{theorem}

Based on Theorem~\ref{thm:Grunbaum}, let $\theta(K) = \sup \left\{ \frac{\vol_d(K_u)}{\vol_d(K_{-u})} : u \in \Sph^{d-1}   \right\}$.
This quantity can be regarded as a measure of symmetry of $K$: it satisfies $\theta(K) \geq 1$ for every convex body $K$, where we have equality for centrally symmetric bodies, and it attains its maximal value $\theta_d = \left( 1 + \frac{1}{d} \right)^{d} -1$ e.g. for convex cones.

In 2000, Groemer \cite{Groemer2000} proved a stability version of Gr\"unbaum's result. Namely, he proved that for every $\varepsilon \geq 0$, there is some $\lambda_d > 0$ depending only on $d$ such that if $\theta(K) \geq \theta_d - \varepsilon$ for some convex body $K$ in $\Re^d$, then there is a convex cone $C$ in $\Re^d$
such that
\[
\vol_d (K \Delta C) \leq  \lambda_d \vol_d(K) \varepsilon^{1/(2d^2)},
\]
where $X \Delta Y$ denotes the symmetric difference of $X, Y$, i.e. the set of those points which are either in $X$ or in $Y$, but not in their intersection.

Using the inequality $\left( 1 + \frac{1}{d} \right)^{d} < e$, Theorem~\ref{thm:Grunbaum} yields a simple dimension-free estimate for the volumes of the intersections: for any convex body $K$ in $\Re^d$, and for any closed half space $H$ containing the centroid of $K$ in its boundary, $\vol_d(H \cap K ) \geq \frac{1}{e} \vol_d(K)$. This inequality was generalized by Bertsimas and Vempala \cite{BV2004} for a convex body $K$ in \emph{isotropic position}, defined as a special type of centered convex body: for any unit vector $u \in \Sph^{d-1}$, $K$ satisfies the equality $1 = \frac{1}{\vol_d(K)} \int_K (\langle u,x \rangle)^2 \, dx$. Their result states that if $K$ is a convex body in $\Re^d$ in isotropic position, and $z \in K$ is at distance $t$ from the origin, then for any closed half space $H$ containing $z$, $\vol_d(H \cap K) \geq \left( \frac{1}{e} - t \right) \vol_d(K)$. Here we note that for any centered convex body $K$ there is a linear transformation $A$ such that $A(K)$ is in isotropic position (see e.g. \cite{BV2004}). Since linear transformations do not change the ratio of volumes and the lengths of parallel segments, it seems possible to restate the result in \cite{BV2004} for an arbitrary convex body $K$ centered at $o$. This was done in \cite{SY2021} by Shyntar and Yaskin, who strengthened of the result in \cite{BV2004} for $\Re^2$ in the following form: Let $K$ be a plane convex body with $o$ as its centroid. Let $-1 < \alpha < 2$ and $u \in \Sph^1$. Let $H$ be the closed half plane defined as $H= \{ x \in \Re^2 : \langle x , u \rangle \geq \alpha h_K(-u) \}$. Then
\begin{equation}\label{eq:SY1}
C_1(\alpha) \area(K) \leq \area(K \cap H) \leq C_2(\alpha) \area(K),
\end{equation}
where
\[
C_1(\alpha) = \left\{
\begin{array}{l}
\frac{1}{9}(2-\alpha)^2, \hbox{ if } \alpha \in (-1,0),\\
\frac{4}{9}(1+\alpha)(1-2\alpha), \hbox{ if } \alpha \in (0,1/2),\\
0, \hbox{ if } \alpha \in (1/2,2),
\end{array}
\right.
\]
and
\[
C_2(\alpha) = \left\{
\begin{array}{l}
1-\frac{4}{9}(1+\alpha)^2, \hbox{ if } \alpha \in (-1,0),\\
\frac{5-3\alpha}{9(1+\alpha)}, \hbox{ if } \alpha \in (0,1),\\
\frac{1}{9}(2-\alpha)^2, \hbox{ if } \alpha \in (1,2).
\end{array}
\right.
\]
The authors of \cite{SY2021} also characterized the equality cases in (\ref{eq:SY1}).

Observe that if $K$ is a centered convex body in $\Re^d$, and $H$ is a hyperplane through $o$, then $o$ may not be the centroid of $K | H$. Thus, it is a meaningful question to examine Gr\"unbaum's problem for projections of convex bodies.
This problem was completely solved by Stephen and Zhang \cite{SZ2017}, who proved that for any centered convex body $K$ in $\Re^d$, any $k$-dimensional linear subspace $E$ with $1 \leq k \leq d$, and any half space $H$ bounded by $u^{\perp}$ with $u \in E \cap \Sph^{d-1}$, we have
\[
\vol_k ((K|E) \cap H) \geq \left( \frac{k}{d+1} \right)^k \vol_k (K | E).
\] 
Furthermore, here there is equality for some $K$ and $E$ if and only if $K = \conv \{ y_1 + L_1, y_2 + L_2\}$, where
\begin{itemize}
\item $L_1 \subset u^{\perp}$ and $L_1 | (E \cap u^{\perp}$ are $(k-1)$-dimensional convex bodies,
\item $L_2 \subset E^{\perp}$ is a $(d-k)$-dimensional convex body,
\item $y_1, y_2 \in \Re^d$ satisfy $y_1 \notin H$ and $y_2 \in \inter(H)$,
\item the centroid of $K$ is $o$.
\end{itemize}
More specifically, the authors of \cite{SZ2017} prove a more general version of the above result for volumes, where the centroid of $K$ is replaced by an arbitrary point in the interior of $K$, representing the center of mass of an inhomogeneous object.

In the same vein, one can examine Gr\"unbaum's problem for sections of a convex body. 
This was proposed by Fradelizi, Meyer and Yaskin \cite{FMY2017}, who proved that for any centered convex body $K$, linear subspace $V$ of dimension $d-k$ with $0 \leq k \leq d-1$, and closed half space $H$ with boundary $u^{\perp}$ for some $u \in V \cap \Sph^{d-1}$, we have
\begin{equation}\label{eq:GRsection}
\vol_{d-k}(K \cap V \cap H) \geq c  \vol_{d-k}(K \cap V),
\end{equation}
where $c \geq \frac{c_0}{(k+1)^2} \left( 1 + \frac{k+1}{d-k} \right)^{-(d-k-2)}$ for some absolute constant $c_0 > 0$.
This result follows from a more general setting investigated in \cite{FMY2017} for intersections of a convex body $K$ with a convex cone, which the authors use to prove a conjecture of Meyer and Reisner \cite{MR2011} about the volumes of convex intersection bodies.
Gr\"unbaum's problem for sections was settled in \cite{MSZ2018} by Myroshnychenko, Stephen and Zhang, who proved the following: for any convex body $K$ in $\Re^d$, a $k$-dimensional linear subspace $E$ of $\Re^d$ and a closed half space $H$ bounded by $u^{\perp}$ with $u \in \Sph^{d-1} \cap E$, if the centroid of $K$ lies in $E \cap u^{\perp}$, then
\[
\vol_k(K \cap E \cap H) \geq \left( \frac{k}{d+1} \right)^k \vol_k (K \cap E).
\]
Furthermore, here we have equality if and only if $K = \conv \left\{ -\left( \frac{d-k+1}{k} \right) z + D_0, z+ D_1 \right\}$, where
\begin{itemize}
\item $z \in E \cap \inter (H)$,
\item $D_0$ is a $(k-1)$-dimensional convex body in $E \cap u^{\perp}$,
\item $D_1$ is a centered $(d-k)$-dimensional convex body in a $(d-k)$-dimensional subspace $F$ of $\Re^d$ such that $E \cup F$ spans $\Re^d$.
\end{itemize}

The results about both the intersections and the projections were generalized for dual volumes by Stephen and Yaskin \cite{SY2019} as follows.
Let $K$ be a centered convex body in $\Re^d$, and let $E$ be a $k$-dimensional linear subspace of $\Re^d$, where $1 \leq k \leq d$. Let $H$ be a closed half space bounded by $u^{\perp}$ with $u \in \Sph^{d-1} \cap E$. Then, for any $1\leq i \leq k$,
\[
\tilde{V}_i ((K | E) \cap H) \geq \left( \frac{i}{d+1} \right)^i \tilde{V}_i(K | E),
\]
and
\[
\tilde{V}_i (K \cap E \cap H) \geq \left( \frac{i}{d+1} \right)^i \tilde{V}_i(K \cap E),
\]
where $\tilde{V}_i(\cdot)$ denotes $i$th dual volume with respect to the $k$-dimensional subspace $E$.
It is worth noting that the case $i=k$ of these statements coincides with the results in \cite{SZ2017} and \cite{MSZ2018}, respectively, mentioned in the previous two paragraphs.

An interesting discrete version of the theorem of Winternitz was proposed by Shyntar and Yaskin \cite{SY2021}. To state the problem, for any set $X\subseteq \Re^2$ we denote the cardinality of the point set $X \cap \mathbb{Z}^2$ by $\# X$, and we call a convex polygon $P$ an \emph{integer polygon} if all the vertices of $P$ have integer coordinates. In this setting, the problem of Winternitz asks for finding the largest positive number $C>0$ such that for any integer polygon $P$, and any closed half plane $H$ containing $c(P)$ in its boundary, we have $\frac{\# (K \cap H)}{\# K} \geq C$. The authors of \cite{SY2021} give an example showing that there is no positive number $C$ that satisfies this property, and prove that, for any integer polygon $P$ and any closed half plane containing $c(3P)$ in its boundary,
\[
\frac{\# ((3P) \cap H)}{\# (3P)} > \frac{1}{6}.
\]
Furthermore, for any $t \geq 6$, any integer polygon $P$ and any closed half space containing $c(tP)$ in its boundary, we have
\begin{equation}\label{eq:discreteWinternitz}
\frac{\# ((tP) \cap H)}{\# (tP)} \geq \frac{\frac{4}{9}t^2-2t-3}{t^2+3t+2}.
\end{equation}
Here we note that for $t >0$, the expression on the right-hand side of (\ref{eq:discreteWinternitz}) is increasing, and tends to $\frac{4}{9}$ as $t \to \infty$.

A generalization of Gr\"unbam's theorem for not necessarily convex sets was considered by Mar\'\i n Sola and Yepes Nicol\'as \cite{MSYN2021}.
To state their result, we recall that by Brunn's theorem, for any compact, convex set $K$ in $\Re^d$, and any $u \in \Sph^{d-1}$, the function $f_X : \Re \to \Re$, $f(t) = \vol_{d-1}(K \cap (tu + u^{\perp}))$ is $\frac{1}{d-1}$-concave on its support. The authors of \cite{MSYN2021} proved that if $X \subset \Re^d$ is a centered compact set with nonempty interior, and $H$ is a closed half space with $u^{\perp}$ as its boundary for some $u \in \Sph^{d-1}$, and the function $f_X : \Re \to \Re$, $f(t) = \vol_{d-1}(X \cap (tu + u^{\perp}))$ is $p$-concave on its support for some $p > 0$, then
\[
\frac{\vol_d(X \cap H)}{\vol_d(X)} \geq \left( \frac{p+1}{2p+1} \right)^{(p+1)/p},
\]
and if $f_X$ is log-concave on its support, then
\[
\frac{\vol_d(X \cap H)}{\vol_d(X)} \geq \frac{1}{e}.
\]
They also characterize the equality case. The same problem, under the more general condition that the function $f_X$ is $\varphi$-concave for some strictly increasing, continuous function $\varphi: [0,\infty ) \to \Re \cup \{ -\infty \}$ satisfying some mild regularity conditions, was investigated by Mar\'\i n Sola \cite{MS2024}.

\subsection{The Busemann-Petty centroid inequality}\label{subsec:BP}

Let $K$ be a centered convex body in $\Re^d$. Consider the function $H_K: \Re^d \to \Re$, defined by
\begin{equation}\label{eq:centroidbody}
H_K(u) = \frac{1}{\vol_d(K)} \int_K | \langle x, u \rangle | \, dx.
\end{equation}
It is well known that this function is positively homogeneous, nonnegative, convex function which is zero only at $o$ \cite{Petty1961}. Thus, $H_K$ is the support function of a convex body $CK$ (see Section~\ref{sec:intro}), called the \emph{centroid body} of $K$.
This body is strictly convex, and if $K$ is symmetric with respect to $o$, i.e. $-K=K$, then the hypersurface $\bd(CK)$, called the \emph{centroid surface} of $K$, coincides with the set of the centroids of the intersections of $K$ with the closed half spaces containing $o$ in their boundaries. Centroid bodies, attributed by Blaschke to Dupin (see e.g. \cite{Schneiderbook}), were given their name in the seminal work of Petty \cite{Petty1961}, who proved the inequality that
\begin{equation}\label{eq:BusemannPetty}
\vol_d(CK) \geq \left( \frac{(d+1) \kappa_d}{2 \kappa_{d-1}} \right)^d \vol_d(K),
\end{equation}
with equality if and only if $K$ is an $o$-symmetric ellipsoid. This inequality is called the \emph{Busemann-Petty centroid inequality}. In 2016, Ivaki \cite{Ivaki2016} proved a stability version of the equality case of this inequality in the plane. Namely, he proved that there are positive constants $\varepsilon_0$ and $\gamma$ such that for any $0 < \varepsilon < \varepsilon_0$ and any plane convex body $K$ satisfying $\frac{\area(CK)}{\area(K)} \leq \left( \frac{4}{3\pi} \right)^2 (1+\varepsilon)$, the Banach-Mazur distance of $K$ and an ellipsoid is at most $1+ \gamma \varepsilon^{1/8}$.

It was observed by Milman and Pajor \cite{MP1989} that all centroid bodies are zonoids; that is, they can be obtained as limits, in the Hausdorff distance, of a sequence of convex bodies, called zonotopes, obtained as the Minkowski sums of finitely many closed segments.

The concept of centroid body, using a different normalization, was extended to $L_p$-Brunn-Minkowski theory by Lutwak and Zhang \cite{LZ1997} in the following way. For any real $p \geq 1$, let
\[
c_{d,p} = \frac{\kappa_{d+p}}{\kappa_2 \kappa_d \kappa_{p-1}}.
\]
For any star body $K \subset \Re^d$, let $\Gamma_p^* K$ be the $o$-symmetric convex body defining the norm
\[
||x||_{\Gamma_p^* K} = \left( \frac{1}{c_{n,p} \vol_d(K)} \int_K | \langle x,y \rangle |^p \, d y \right)^{1/p}.
\]
This body is called the \emph{polar $L_p$-centroid body} of $K$, and it can be extended to the case $p=\infty$ by taking the limit $p \to \infty$. It is worth noting that if $K$ is an $o$-symmetric convex body and $p=\infty$, then $\Gamma_p^* K$ coincides with the Euclidean polar body $K^*$ of $K$, defined as $K^* = \{ y : \langle x,y \rangle \leq 1 \hbox{ for every } x \in K \}$. It was proved in \cite{LZ1997} that for any star body $K$, we have
\begin{equation}\label{eq:LZ1997}
\vol_d(K) \vol_d(\Gamma_p^* K) \leq \kappa_d^2,
\end{equation}
with equality if and only if $K$ is an $o$-symmetric ellipsoid. This inequality, for the case $p=\infty$, coincides with the Blaschke-Santal\'o inequality.
The polar body of $\Gamma_p^* K$ is called the \emph{$L_p$-centroid body of $K$}; if $p=1$, this body coincides with a normalization of the centroid body $CK$ of $K$. We note that the body $\Gamma_2 K$ also has a geometric meaning: this is the ellipsoid which has the same inertia about every axis as $K$. This body is called the \emph{ellipsoid of inertia} or \emph{Legendre ellipsoid} of $K$ \cite{CG2002}. Verifying a conjecture in \cite{LZ1997}, Lutwak, Yang and Zhang \cite{LYZ2000} proved that
\begin{equation}\label{eq:LYZ2000}
\vol_d(K) \leq \vol_d(\Gamma_p K),
\end{equation}
with equality if and only if $K$ is an $o$-symmetric ellipsoid. Note that (\ref{eq:LYZ2000}), combined with the Blaschke-Santal\'o inequality, yields (\ref{eq:LZ1997}), and that (\ref{eq:LYZ2000}) is a generalization of the Busemann-Petty centroid inequality in (\ref{eq:BusemannPetty}). A different proof of the inequality in (\ref{eq:LYZ2000}) was given in \cite{CG2002}, and a version, establishing a stronger inequality for nonsymmetric bodies, can be found in the paper \cite{HS2009} of Haberl and Schuster. The estimate in \cite{HS2009} was generalized for compact sets in \cite{QB2014}. A different stronger version of the $L_p$ Busemann-Petty centroid inequality can be found in \cite{MZ2015}, and a variant for the first quermassintegrals of the bodies in \cite{Wang2022}.

Let $\phi : \Re \to [0,\infty )$ be a convex function satisfying $\phi(0)=0$ such that $\phi$ is strictly decreasing on $(-\infty,0]$ or strictly increasing on $[0,\infty)$. The \emph{Orlicz centroid body} $\Gamma_{\phi} K$ of a star body $K$ is defined as the convex body with support function
\begin{equation}\label{eq:Orliczcentroid}
h_{\Gamma_{\phi} K}(x) = \inf \left\{ \lambda > 0 : \frac{1}{\vol_d(K)} \int_K \phi \left( \frac{ \langle x,y \rangle }{\lambda}\right) \leq 1 \right\}.
\end{equation}
We note that if $\phi(t) = |t|^p$ with $p \geq 1$, then $\Gamma_{\phi} K = \Gamma_p K$.
This concept was introduced by Lutwak, Yang and Zhang in \cite{LYZ2010}, who proved that for any function $\phi$ satisfying the above conditions, the ratio $\frac{\vol_d(\Gamma_{\phi} K)}{\vol_d(K)}$ is minimal on the family of convex bodies $K$ with $o \in \inter (K)$ if and only if $K$ is an $o$-symmetric ellipsoid.
This result was extended for arbitrary star bodies by Zhu \cite{Zhu2012}, and the equality case, under a weaker condition on $\phi$, was studied in \cite{WZ2018}.

It is an interesting question to find a converse inequality to (\ref{eq:BusemannPetty}).
This question was proposed by Bisztriczky and B\"or\"oczky Jr. \cite{BB2001} in 2001. They proved that in the family of $o$-symmetric plane convex bodies (resp. plane convex bodies containing $o$), the quantity $\frac{\area(CK)}{\area(K)}$ is maximal if $K$ is a parallelogram (resp. a triangle with $o$ as a vertex). Campi and Gronchi \cite{CG2002_2} extended these results to $L_p$-centroid bodies in the plane, and Chen, Zhou and Yang \cite{CZY2011} to Orlicz centroid bodies in the plane.

The definition of Orlicz centroid bodies in (\ref{eq:Orliczcentroid}) can be extended to that of Orlicz-Lorenz centroid bodies by adding a weight function to the integral in (\ref{eq:Orliczcentroid}). The Busemann-Petty centroid inequality in this setting was investigated in \cite{Nguyen2018, CY2019, FM2019, ZF2019, Zhang2023, FLM2023}. Variants of centroid bodies, and the Busemann-Petty centroid inequality, in non-Euclidean spaces (including complex, spherical and hyperbolic spaces) can be found in \cite{Haberl2019, LW2020, JL2021, BHPS2023, LZS2023}.


\section{Static equilibrium points}\label{sec:equilibria}

Consider a convex body $K$ in $\Re^3$. Clearly, if $K$ is balanced on a horizontal plane $H$, then it is in static equilibrium if and only if the orthogonal projection of its center of mass onto $H$ belongs to $K$. This leads to the following definition (see e.g. \cite{Langi_2022}).

\begin{definition}\label{defn:equilibrium}
Let $K$ be a convex body in $\Re^d$ and $c \in \Re^d$. A point $q \in \bd (K)$ is called an \emph{equilibrium point} of $K$ \emph{with respect to $c$} if the hyperplane passing through $q$ and perpendicular to $[c,q]$ supports $K$. 
\end{definition}

We note that by convention 
every point $c \in \bd(K)$ is an equilibrium point of $K$ with respect to $c$. In the literature, $c$ is often assumed to be the centroid of $K$. In this case we call $q$ simply an \emph{equilibrium point} of $K$. As we mentioned in the beginning of Section~\ref{sec:centroid}, any reference point $c \in \inter(K)$ is the center of mass of $K$ equipped with some suitable, possibly inhomogeneous density function.
Thus, unless we state it otherwise, we assume from now on that the reference point lies in the interior of the body.

For the case $c=o \in \inter(K)$, a simple description of equilibrium points can be given by means of the radial or the support function of $K$. More specifically, if $\rho_K$ or $h_K$ is $C^1$-class, then $q$ is an equilibrium point of $K$ with respect to $o$ if and only if $\frac{q}{|q|}$ is a critical point of $\rho_K$ or $h_K$, respectively. 
Indeed, since the restriction of $\rho_K$ to a plane coincides with the radial function of the intersection of $K$ with the plane, and a similar observation holds for $h_K$ and the projection of $K$ onto the plane, it is sufficient to show our remark for plane convex bodies. For this case a simple geometric argument shows the statement (see also \cite{domokos2016topological, allemann2021equilibria}).
 We note that the properties that $h_K$ or $\rho_K$ are $C^1$-class are not equivalent \cite{Schneiderbook}, but if both are $C^1$-class, their critical points coincide.

A remarkable result of Zamfirescu \cite{zamfirescu1995convex} states that, in Baire category sense with respect to the topology induced by Hausdorff distance, a typical convex body has infinitely many equilibrium points with respect to a typical point. On the other hand, Montejano \cite{montejano1974problem} showed that the Euclidean balls centered at $o$ are the only convex bodies whose every boundary point is an equilibrium point with respect to $o$.

Next, we define non-degenerate equilibria.
\begin{definition}\label{defn:nondegenerates}
Let $K$ be a smooth convex body in $\Re^d$ and $c\in \Re^d$. Assume that $K$ has finitely many equilibrium points with respect to $c$, and its boundary is $C^2$-class in the neighborhood of each equilibrium point $q \in \bd (K)$. Assume that the Hessian of the function $\delta_K : \bd(K) \to \Re$, $\delta_K(x)=|x-c|$ is not zero at any equilibrium point $q$. Then we say that $K$ is \emph{nondegenerate} with respect to $c$. Furthermore, the \emph{index} of any equilibrium point is the number of negative eigenvalues of the Hessian of $\delta_K$ at $q$. Equilibrium points of index $0$ and $d-1$ are called \emph{stable} and \emph{unstable} points, respectively, and equilibrium points of all other indices are called \emph{saddle} points.
\end{definition}

We note that for the special case $c=o\in\inter(K)$, there is an equivalent formulation of the concepts in Definition~\ref{defn:nondegenerates} in terms of the radial function $\rho_K$. Indeed, the conditions for the nondegeneracy of $K$ are equivalent to the properties that $\rho_K$ is $C^1$-class, it has finitely many critical points and $\rho_K$ is $C^2$-class in a neighborhood of every one of them, and at every critical point the Hessian of $\rho_K$ is not zero. Furthermore, the index of the equilibrium point $q \in \bd(K)$ is equal to the number of the negative eigenvalues of the Hessian of $\rho_K$ at the corresponding critical point $\frac{q}{|q|}$.

Whereas there is no definition, in general, for a not necessarily smooth convex body being nondegenerate, the condition in Definition~\ref{defn:nondegenerates} can be replaced by a geometric one if $K$ is a polytope. We note that by definition, every face (of any dimension) of a convex polytope contains at most one point of equilibrium in its relative interior.

\begin{definition}\label{defn:nondegeneratep}
Let $P$ be a convex polytope in $\Re^d$, and let $c \in \Re^d$. For any equilibrium point $q$ of $P$ with respect to $c$, let $H_q$ denote the supporting hyperlane of $P$ at $q$ perpendicular to $[c,q]$. Set $F_q= P \cap H_q$, and let $f_q = \dim F_q$. If for every equilibrium point $q$, $q$ lies in the relative interior of $F_q$, we say that $P$ is \emph{nondegenerate}. In this case the \emph{index} of the equilibrium point $q$ is the quantity $(d-1)-f_q$.
\end{definition}


Consider a $C^1$-class function $f : \Sph^{d-1} \to \Re$ with finitely many critical points such that $f$ is $C^2$-class in a neighborhood of every critical point, and its Hessian is nonzero at every critical point. A standard convolution technique (see e.g. \cite{Ghomi2002}) shows that there is a $C^{\infty}$-class function $g : \Sph^{d-1} \to \Re$ arbitrarily close to $f$ such that for every $0 \leq k \leq d-1$, the number of critical points of $f$ with index $k$ coincide with the same quantity for $g$.
The function $g$ is a Morse function, and hence, applying this observation for the radial function $\rho_K$ of a nondegenerate smooth convex body $K$ with $c=o \in \inter (K)$, and combining it with the Poincar\'e-Hopf formula for the sphere $\Sph^{d-1}$, we have the following.

\begin{theorem}\label{thm:PoincareHopf}
Let $K$ be a nondegenerate smooth convex body in $\Re^3$ with respect to $c \in \inter (K)$, and let $S, H, U$ denote the numbers of stable, saddle and unstable points of $K$ with respect to $c$. Then
\begin{equation}\label{eq:PoincareHopf}
S-H+U=2.
\end{equation}
Similarly, if $K$ is a nondegenerate plane convex body with respect to some $c \in \inter (K)$, with $S$ stable and $U$ unstable points, then $S=U$.
\end{theorem}

Finally, let $P$ be a nondegenerate convex polytope in $\Re^d$ with respect to $c \in \inter (P)$. Then, if $\varepsilon > 0$ is sufficiently small, the convex body $K=P + \varepsilon \BB^d$ is a nondegenerate convex body with respect to $c \in \inter(K)$, and for every $0 \leq k \leq d-1$, the number of equilibrium points of $P$ of index $k$ coincide with the same quantity for $K$. Thus, Theorem~\ref{thm:PoincareHopf} is satisfied also for nondegenerate convex polyhedra and polygons. We also note that the notion of radial function and Theorem~\ref{thm:PoincareHopf} can be extended to the case $c \notin \inter(K)$ in a natural way \cite{gardner1995geometric}.

\subsection{Primary classes of convex bodies based on their equilibria}\label{subsec:primary}

\begin{figure}
    \centering
    \includegraphics[width=\textwidth]{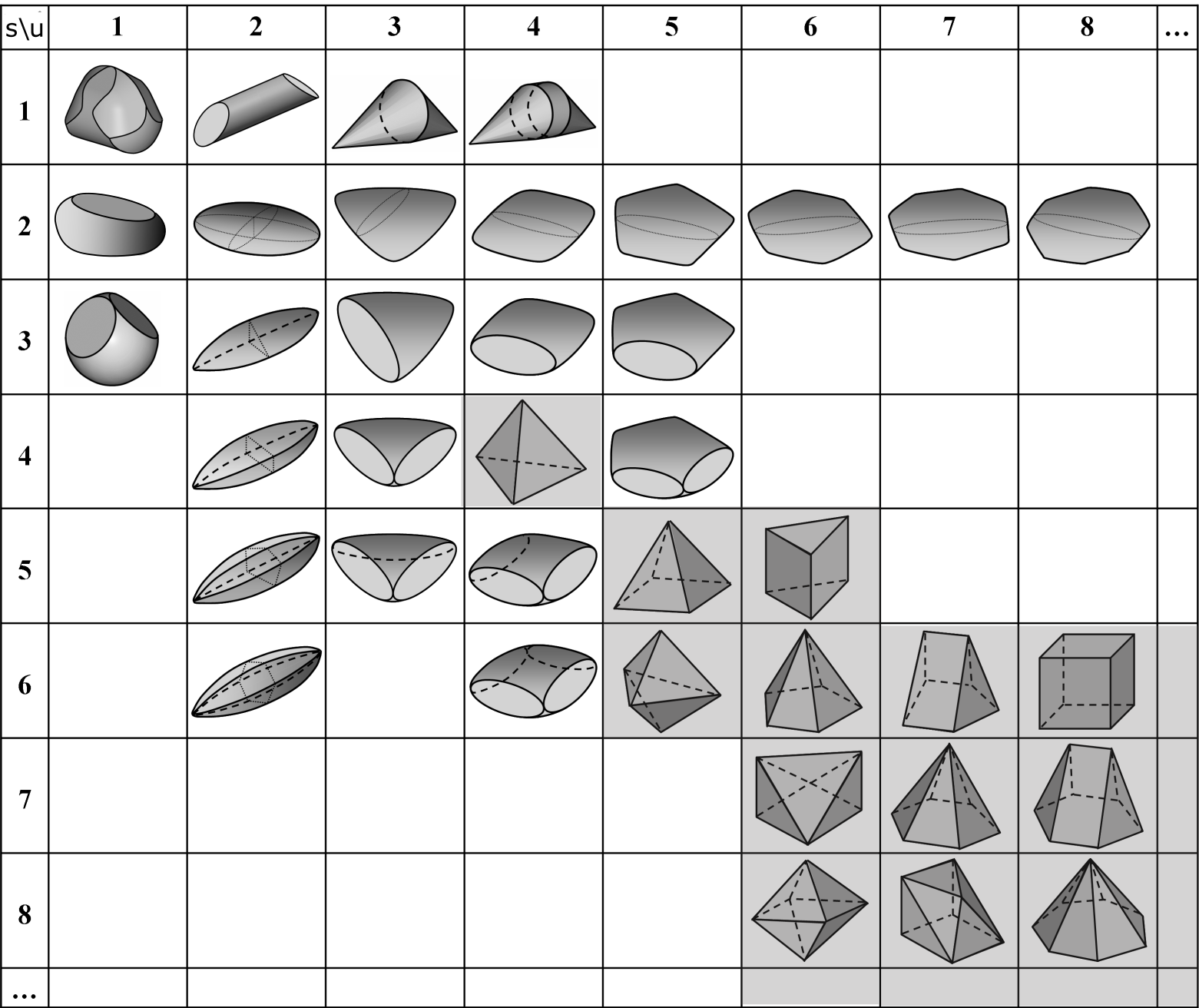}
    \caption{Equilibrium classes with some representative examples. Classes corresponding to polyhedral pairs (see Section \ref{subsec:polyhedra}, below) are highlighted by grey background. Unistable, uni-unstable, and mono-monostatic objects correspond to the first row, the first column, and the top-left corner of the table. The examples have been provided by T\'imea Szab\'o and G\'abor Domokos.}
    \label{fig:table}
\end{figure}

One can classify convex bodies based on their static equilibrium properties by using the numbers of their equilibrium points of different indices. We observe that by Theorem~\ref{thm:PoincareHopf}, for any nondegenerate plane convex body it is sufficient to know the number of either the stable or the unstable points, and in $3$-dimensions two of the numbers of stable, saddle and unstable points. Thus, we define the set $\{ S  \}$ as the family of all nondegenerate convex bodies in $\Re^2$ that have $S$ stable points, and 
the set $\{ S, U \}$ as the family of all nondegenerate convex bodies in $\Re^3$ that have $S$ stable and $U$ unstable points. For later use, we decompose each class $\{ S \}$ into two subclasses $\{ S \}_s$, $\{ S \}_p$, consisting of the nondegenerate smooth plane convex bodies and the convex polygons in this class, respectively, and we introduce the classes $\{ S, U \}_s$, $\{ S, U \}_p$ in the same way (Figure \ref{fig:table}).
Elements of class $\{ 1, U \}$ for some $U$ are called \emph{monostable} or \emph{unistable}, elements of class $\{ S, 1 \}$ are called \emph{mono-unstable} or \emph{uni-unstable}, and a body that is both monostable and mono-unstable is called \emph{mono-monostatic}.

We call this classification the primary classification system, and a class $\{ S \}$ or $\{ S, U \}$ the \emph{primary class} of its elements \cite{varkonyi2006static, domokos2016topological}. It is a natural question to ask what are the integers $S \geq 0$ (resp. the integers $S,U \geq 0$) such that the class $\{ S \}$ (resp. the class $\{ S,U \}$) is not empty.
Our first step in answering this question is the observation that if $K$ is a nondegenerate convex body, then every boundary point of $K$ closest to $c(K)$ is a  stable, and every boundary point of $K$ farthest from $c(K)$ is an unstable point of $K$. This implies that in every nonempty primary class, $S,U $ are positive.
The first nontrivial result in this direction is due to Domokos, Papadopulos and Ruina \cite{domokos1994static}, who proved that the class $\{1\}$ is empty, or in other words, every nondegenerate plane convex body has at least two stable and two unstable points. On the other hand, it is easy to see that the class $\{ S \}$ is not empty for any $S \geq 2$. Indeed, the class $\{ 2 \}$ contains all ellipses different from a circle, and for any $S \geq 3$, the class $\{ S \}$ contains, e.g., the regular $S$-gons. This completes the planar case. As a final remark, we mention that the result that $\{ 1 \}$ is empty is equivalent to the famous $4$-vertex theorem in differential geometry \cite{varkonyi2006static}.

\begin{figure}
    \centering
    \includegraphics[width=1\linewidth]{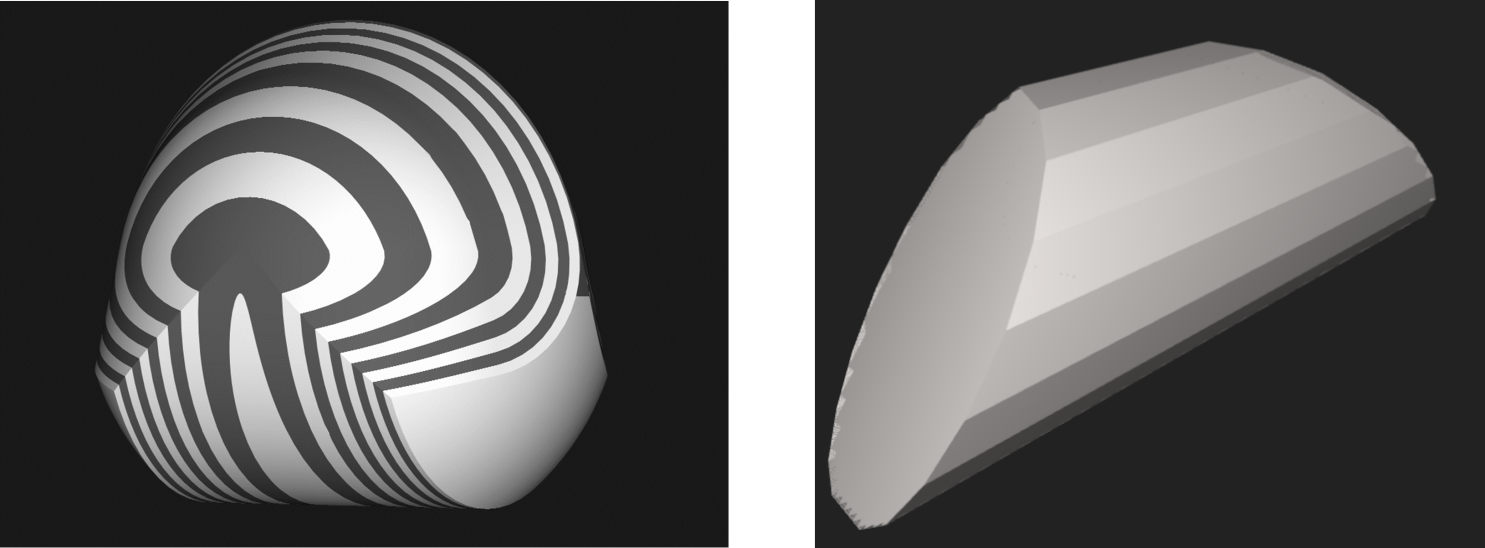}
    \caption{Left: illustration of the G\"omb\"oc: a piecewise smooth body in class $\{1,1\}$. The boundary is coloured according to level curves of the radial function, which has two critical points at the top and at the bottom of the surface. Right: the monostatic polyhedron of Guy \cite{conway1966stability}, see Section \ref{subsec:polyhedra}, below. The illustration has been created with the aid of a 3D model available at \emph{https://www.thingiverse.com/thing:90866}.
    \label{fig:GombocGuy}
    }
\end{figure}

The characterization of nonempty primary classes $\{ S, U \}$ was carried out by  V\'arkonyi and Domokos~\cite{varkonyi2006static}. In particular, they presented a construction of a convex body with only one stable and one unstable point. For preciseness, we mention that the body constructed in \cite{varkonyi2006static}  was not $C^2$-class at its two equilibrium points, and thus, it is not nondegenerate as defined in Definition~\ref{defn:nondegenerates}; nevertheless the radial function of the body has a global minimum and maximum at its equilibrium points. We also add that, as it was remarked in \cite{domokos2016topological}, a standard smoothing technique yields a $C^{\infty}$-class, nondegenerate approximation of this body with the same numbers of stable and unstable points. For an explicit construction of nondegenerate convex bodies in class $\{ 1,1 \}$, with $C^2$-class boundaries, the reader is also referred to \cite{gabor2023characterization}. A representative example of class $\{1,1\}$ constructed by patching together elementary surfaces became known as `G\"omb\"oc' (Figure \ref{fig:GombocGuy}). The existence of a convex polyhedron in class $\{ 1,1 \}$ was shown in \cite{Langi_2022}. Another remarkable simple shape in class $\{ 1,1 \}$ given by an explicit formula was recently published by Sloan \cite{sloan2023analytical}, and possible symmetries of mono-monostatic bodies have been determined in \cite{gabor2023characterization}.

Once it has been verified that the class $\{ 1,1 \}$ is not empty, one can also show that no class $\{ S, U \}$, with $S,U \geq 1$ is empty by finding a way to increase the number of equilibrium points of a body. This was also carried out in \cite{varkonyi2006static} by Domokos and V\'arkonyi, who showed that, in a neighborhood of a fixed stable point of a convex body $K$, it is possible to create a stable point by truncating $K$ by a suitable plane, or an unstable point by truncating $K$ by a suitable cone. These steps are called `Columbus steps' in memory of the famous discoverer Christopher Columbus, who balanced an egg on its top by creating a stable point near it. Since the class $\{ 1,1 \}$ is not empty, subsequent applications of these steps show that no class $\{ S, U \}$ with $S,U \geq 1$ is empty. We note that a more rigorous investigation of Columbus steps, involving also a construction of a smooth approximation of the resulting body, can be found in \cite{domokos2016topological}. The fact that all these classes contain also nondegenerate convex polyhedra follows from \cite{Langi_2022}. 

\subsection{Secondary and tertiary classes}\label{subsec:sectert}

In this subsection, we investigate nondegenerate smooth convex bodies in $\Re^3$. Consider the radial function $\rho_K : \Sph^2 \to \Re$ of such a body. 
By the smoothing technique described in the beginning of Section~\ref{sec:equilibria}, we may assume that $\rho_K$ is a $C^{\infty}$-class function. Then $\rho_K$ is a Morse function on $\Sph^2$, and the gradient of $\rho_K$ is nonzero at $q$ if and only if $q$ is not critical, i.e. it is regular. By the Picard-Lindel\"of theorem, for any regular point $q$ of $\Sph^2$ there is a unique curve $\gamma: \Re \to \Sph^2$ with the property that $\gamma(0)=q$, and $\gamma'(t)$ is equal to the gradient of $\rho_K$ at $\gamma(t)$ for all $t \in \Re$. Such a curve $\gamma$, called an \emph{integral curve} of the gradient flow induced by $\rho_K$ starts and ends at two different critical points of $\rho_K$, called the \emph{origin} and the \emph{destination} of $\gamma$. For any critical point $x$, the \emph{descending (resp. ascending) manifold} of $x$, denoted by $D(x)$ (resp. $A(x)$) is the union of $\{x\}$ and all integral curves with $x$ as their destination (resp. origin). The connected components of the nonempty sets, obtained as the intersection of an ascending and a descending manifold, are called \emph{cells} of the \emph{Morse-Smale complex generated by $\rho_K$}. The $2$-, $1$-, and $0$-dimensional cells of this complex are called \emph{faces, edges}, and \emph{vertices}, respectively.

If $x \neq y$ are critical points of $\rho_K$, then $A(x) \cap D(x) = \{ x \}$, and $A(x) \cap D(y)$ is the union of all integral curves with origin $x$ and destination $y$. Thus, the vertices of the Morse-Smale complex of $\rho_K$ are the critical points of $\rho_K$, and if $x$ is a stable and $y$ is an unstable point, then $A(x) \cap D(y)$ is an open set in $\Sph^2$, implying that its connected components are faces of the complex (provided ${A(x)\cap D(y)\ne\emptyset}$). On the other hand, the same does not hold if $x$ or $y$ is a saddle point since every saddle point is the origin, and also the destination, of exactly two, \emph{isolated} integral curves, corresponding to edges of the complex \cite{edelsbrunner2003hierarchical, LLD2024}.

The function $\rho_K$ is called \emph{Morse-Smale} if all ascending and descending manifolds of $\rho_K$ intersect only transversally; or equivalently, if any pair of intersecting ascending and descending $1$-dimensional manifolds cross. In this case the crossing point of these $1$-manifolds is necessarily a saddle, since crossing at a regular point would contradict the property that any regular point belongs to exactly one integral curve. In the following we only deal with convex bodies whose radial function $\rho_K$ is Morse-Smale. In this case every face of the Morse-Smale complex is bounded by a union of lower dimensional cells, either two disjoint vertices or a cycle of vertices and edges.

The existence of a face whose boundary is a pair of vertices implies that the body $K$ has no saddle point, that is, its primary class is $\{ 1,1 \}$. In this case the Morse-Smale complex is made up of a unique face corresponding to the single connected component of $A(x)\cap D(y)$, where $x$ is the stable and $y$ is the unstable point of $K$, and its boundary consisting of the two critical points $x$ and $y$. In addition, the function $\rho_K$ has no saddle points, no isolated integral curves, therefore its Morse-Smale complex has no edges. In the opposite case the boundary of every face consists of a cycle of connected edges and vertices, and it is known (see e.g. \cite{edelsbrunner2003hierarchical}) that it consists of four edges, and the four critical points in the boundary are a stable, a saddle, an unstable and a saddle point, in this cyclic order. The boundary is possibly glued to itself along vertices and edges.

The topological graph $G$ on $\Sph^2$, whose vertices are the critical points of $\rho_K$, and whose edges are the edges of the Morse-Smale complex, is called the \emph{Morse-Smale graph} generated by $\rho_K$ \cite{domokos2016topological}. This graph is usually regarded as a \emph{$3$-colored} quadrangulation of $\Sph^2$, where the `colors' of the vertices are the three types of a critical point. We regard the topological graph $G$ defined by $\rho_K$ as the drawing of a $3$-colored abstract graph $\bar{G}$ on $\Sph^2$. The graphs $\bar{G}$ and $G$ are called the \emph{secondary} and the \emph{tertiary class} of $K$, respectively. The secondary classification of convex bodies is a refinement of the primary one, as convex bodies whose Morse-Smale complexes induce the same abstract graph have the same numbers of stable, saddle and unstable points, whereas the statement is not necessarily true in the opposite direction. Similarly, tertiary classification is a refinement of the secondary one.

\begin{figure}[ht]
\includegraphics[width=\textwidth]{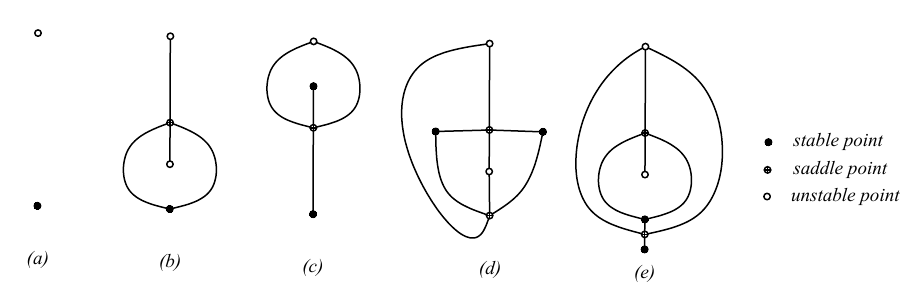}
\caption[]{The topological graphs of convex bodies with at most two stable and two unstable points.}
\label{fig:small_graphs}
\end{figure}

Let us denote by $\Q_3^*$ the family of finite $3$-colored quadrangulations of $\Sph^2$, up to a homeomorphism, where the color classes are denoted by $\mathcal{S}, \mathcal{H}, \mathcal{U}$, each element of $\mathcal{H}$ has degree $4$, and the cardinalities of the classes satisfy the equation $\card \mathcal{S} - \card \mathcal{H} + \card \mathcal{U} = 2$.
We remark that the Morse-Smale complex of every Morse-Smale function on $\Sph^2$ belongs to $\Q_3^*$. Furthermore, let us denote the Morse-Smale graph of a nondegenerate convex body $K$, if it exists, by $Q_3^*(K)$. We note that by convention, the Morse-Smale graphs of the elements of the class $\{ 1,1 \}$, consisting of two vertices with different colors and no edges, are regarded as elements of $\Q_3^*$. These graphs are called \emph{path graphs} on two vertices, and are denoted by $P_2$ (see Figure~\ref{fig:small_graphs}).

As in the case of primary classes, it is a natural question to ask which elements of $\Q_3^*$ appear as Morse-Smale graphs of nondegenerate convex bodies. 
The answer, namely that \emph{every} such quadrangulation is realizable by a convex body, was given in \cite{domokos2016topological} by Domokos, L\'angi and Szab\'o.
The idea of the proof is a refinement of Columbus steps, described in the previous subsection. More specifically, a generalization of a method of Batagelj \cite{Batagelj}, and also Negami and Nakamoto \cite{Negami} shows that every $3$-colored quadrangulation in $\Q_3^*$ can be reduced to $P_2$ by finitely many applications of a graph operation, called \emph{face contraction}, or equivalently, every such quadrangulation can be obtained from $P_2$ by finitely many subsequent \emph{vertex splittings}, the dual operations of face contractions. Based on this result, the authors of \cite{domokos2016topological} show that for any convex body $K$ and any $3$-colored quadrangulation $Q$ in $\Q_3^*$, obtained from $Q_3^*(K)$ by an arbitrary vertex splitting, there is a convex body $K'$, obtained as a smoothened truncation of $K$ by a plane or a cone, satisfying that $Q_3^*(K') = Q$.

\begin{figure}[ht]
\includegraphics[width=\textwidth]{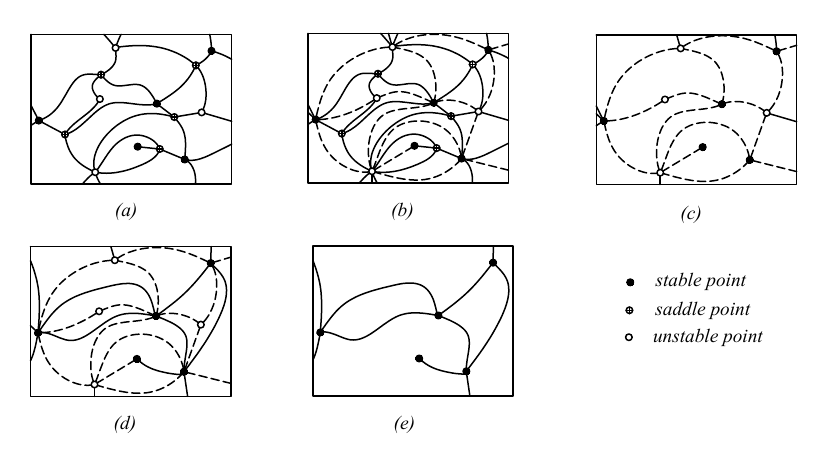}
\caption[]{Different representations of a Morse-Smale complex. a) Primal Morse-Smale graph. b) Primal Morse-Smale graph with the stable and unstable points connected. c) Quasi-dual Morse-Smale graph. d) Quasi-dual graph with the stable points connected. e) The induced uncolored planar graph; vertices, edges and faces correspond to the stable, saddle and unstable points of the primal representation.}
\label{fig:quasi-dual}
\end{figure}

We note that the above representation of the Morse-Smale complex of $K$ is also called the \emph{primal} Morse-Smale graph of $K$.
Saddle points can be removed from the primal Morse-Smale graph without losing information:
first we connect maxima and minima in the quadrangles, then cancel saddle points and edges incident to them (see Figure~\ref{fig:quasi-dual}).
This representation is called the \emph{quasi-dual} Morse-Smale graph (see. \cite{Dong}). Every face of a quasi-dual representation is a quadrangle, and the vertices on its boundary form an alternating sequence of two stable and two unstable points. Connecting the two stable points in each face, and removing all unstable points and their edges we obtain an even more compact representation of the Morse-Smale graph. This representation is a planar graph that might contain loops and multiple edges. The vertices of this graph are the stable points of the Morse-Smale complex, and its edges correspond to the saddle points of the complex (See Figure~\ref{fig:quasi-dual}).

It is a natural question to ask for the numbers of secondary and tertiary classes in each primary class. However, the answer to this question is known only for some small values of $S$ and $U$.


The definitions of the secondary and tertiary classes of a smooth nondegenerate convex body $K$ rely on the topological properties of the Morse-Smale complex of the gradient flow of the radial function $\rho_K$ of the body. This task seems much more complicated if $K$ is a nondegenerate convex polyhedron. Partial results in this direction can be found in the recent paper \cite{LLD2024} of Ludm\'any, L\'angi and Domokos.

\subsection{Transition between different classes}\label{subsec:transition}

The next topic in our survey is to examine how difficult it is to move a convex body from its class to another one.
We start the investigation with primary classes. Following the paper \cite{domokos2014robustness} of Domokos and L\'angi, the question we ask is this: For a given nondegenerate convex body $K$, what is the minimum fraction of the volume of $K$ that is necessary to remove from $K$ to obtain a nondegenerate convex body with a different number of equilibrium points? As we have already seen in Subsection~\ref{subsec:primary}, this number is equal to zero for any smooth convex body, because it is possible to use an arbitrarily small truncation of $K$ to \emph{increase} the number of equilibrium points of the body. This shows that we should aim at \emph{decreasing} the number of equilibrium points of $K$. Thus, we define the \emph{(downward) robustness} of the convex body $K$ in $\Re^d$ as the quantity
\[
\rho(K) = \frac{1}{\vol_d(K)} \inf \left \{ \vol_d (K \setminus K') : K' \in \F_<(K) \right\},
\]
where the set $\F_<(K)$ denotes the family of convex bodies in $K' \subset K$ with fewer equilibrium points than $K$. 
Furthermore, we define the robustness of an equilibrium class as the supremum of the robustnesses of the elements of the class.
To find the robustness of a convex body, or an equilibrium class, is challenging partly because truncating a body changes also the position of its centroid.
This is the motivation behind defining the \emph{internal} and \emph{external robustness} of a convex body. 
As in \cite{domokos2014robustness}, we do it for plane convex bodies with piecewise smooth boundaries, and remark that the nondegeneracy of such a body, and its stable and unstable points, can be defined for them by combining the conditions in Definitions~\ref{defn:nondegenerates} and \ref{defn:nondegeneratep}. We denote the family containing such plane convex bodies by $\K_2$.

\begin{definition}\label{defn:internal}
Let $K \in \K_2$ and $c \in \inter (K)$. Assume that $K$ has $S$ stable points with respect to $c$.
Let $R(K,c) \subseteq \inter (K)$ denote the set of the points such that $K$ has $S$ stable
points with respect to any point of $R(K,c)$.
The \emph{internal robustness} of $K$ with respect to $c$ is
\[
\rin(K,c) = \frac{\inf \left\{ |q-c| : q \notin R(K,c) \right\}}{\perim (K)}  ,
\]
where $\perim (K)$ is the perimeter of $K$.
\end{definition}

The existence of the quantity $\rin(K,c)$ follows from the observations that, first, for any $S > 0$, the set with respect to which $K$ has $S$ stable points is open, and second, Blaschke's rolling ball theorem implies that any $K \in \K_2$ contains two points with respect to which $K$ has different numbers of stable points. The question of how the number of equilibria depends on the location of the center of mass is discussed in detail in Section \ref{subsec:inhom}, and the concept of internal robustness is illustrated by Figure \ref{fig:evolute}.

\begin{definition}\label{defn:external}
Let $K \in \K_2$ and $c \in \inter (K)$. Assume that $K$ has $S$ stable points with respect to $c$. We define the
\emph{(downward) external robustness of $K$} with respect to $c$ as the quantity
\[ \rex(K,c) =
\frac{\inf \{ \area(K \setminus K') : K' \in \F_{<}(K,c) \}}{\area(K)},
\]
where $\F_{<}(K,c)$ denotes the family of the nondegenerate plane convex bodies $K' \in \K_2$, satisfying $c \in \inter (K') \subset K$ , which has strictly less than $S$ stable points with respect to $c$. If the set in the numerator is empty, we let $\rex(K,c)=1$.
\end{definition}

The internal and the external robustness of the equilibrium class $\{ S \}$ is defined as the quantity $\rin_S = \sup \{ \rin(K,c(K)) : K \in \{S \}\}$ and 
$\rex_S = \sup \{ \rex(K,c(K)) : K \in \{S \}\}$, respectively. Regarding these quantities, the authors of \cite{domokos2014robustness} prove that for any $S \geq 3$, both $\rin(K,p)$ and $\rex(K,p)$ are maximal over all $K \in \K_2$ and $c \in \inter (K)$ if $K$ is a regular $S$-gon, and $c=c(K)$.

In \cite{domokos2014robustness}, the authors extend the concepts of external and internal robustness to nondegenerate smooth convex bodies and nondegenerate convex polyhedra in $\Re^3$ in the straightforward way, where, in case of internal robustness, they replace the normalizing factor $\perim (K)$ by $\sqrt{\surf(K)}$, where $\surf(\cdot)$ denotes the surface area. They show that any platonic solid in class $\{ S,U \}$ has maximal internal robustness among the elements of $(S,U)$ having the same numbers of faces, edges and vertices, respectively. They also show that the (full) robustness of the classes $\{ 1,1 \}$, $\{ 1,2 \}$ and $\{ 2, 1 \}$ is equal to $1$.

Transitions between equilibrium classes can be described also from the point of view of dynamical systems. The equivalence classes of nondegenerate, smooth convex bodies are generic in the sense that a  sufficiently small $C^2$-class deformation of the body does not change any of its classes. In particular, if $\{  K(t) : t \in [0,1] \}$ is a $C^2$-class differentiable $1$-parameter family of convex bodies, where $K(0)$ and $K(1)$ belong to different equilibrium classes, then there is some 
value $0 < t_0 < 1$, $K(t_0)$, where $K(t_0)$ is degenerate. If, in addition, for any $0 \leq t < t_0$ (resp. $t_0 < t \leq 1$) $K(t)$ belongs to the same equilibrium class as $K(0)$ (resp. $K(1)$), then there is a bifurcation of the corresponding gradient vector fields induced by $\rho_{K(t)}$. In a `generic' case, this bifurcation has codimension $1$ \cite{Holmesbook}. The codimension $1$ bifurcations of gradient vector fields on $\Sph^2$ have two types, called \emph{saddle-node} bifurcations and \emph{saddle-saddle connections}.

In terms of topological graphs, saddle-node bifurcations correspond to vertex splittings/face contractions, while saddle-saddle connections correspond to another type of graph transformations on $3$-colored quadrangulations on $\Sph^2$, called \emph{diagonal slides}. These graph transformations establish a combinatorial connection between primary/secondary or tertiary equilibrium classes, namely we can say that the transition between two classes $\mathcal{C}_0$ and $\mathcal{C}_1$ induced by such a graph transformation $T$ is \emph{combinatorially realizable} if there are $3$-colored quadrangulations $Q_0$ and $Q_1$, belonging to the classes $\mathcal{C}_0$ and $\mathcal{C}_1$, respectively, such that $Q_1$ is obtained from $Q_0$ by $T$.
In other words, we regard the $3$-colored quadrangulations of $\Sph^2$, up to a homeomorphism, as vertices of a `metagraph' $\mathcal{G}$, and connect two vertices of $\mathcal{G}$ by an edge if one of them is induced from the other one by a vertex splitting or a diagonal slide. We call an edge \emph{primary, secondary} or \emph{tertiary} if the two quadrangulations are induced by convex bodies in different primary, secondary or tertiary classes, respectively.

It is a meaningful question to ask which edges of $\mathcal{G}$ are realizable by a generic $1$-parameter family of convex bodies. This problem was partially solved in \cite{domokos2016genealogy}, where it was proved that every primary and secondary edge of $\mathcal{G}$ is realizable. The case of tertiary edges is still open.

\subsection{The location and clustering of equilibria}\label{subsec:location}

Flattish objects like a coin tend to have two stable equilibrium points, and elongated ones like a needle tend to have two unstable ones. This simple observation can be turned into a more precise statement in various ways. For example, it was shown in \cite{DL2018} that for any convex body $K$ in $\Re^d$ with $C^{\infty}$-class boundary and strictly positive Gaussian curvature, and for any orthogonal affinity $f_{\lambda} : \Re^d \to \Re^d$ with affinity ratio $\lambda > 0$, if $\lambda$ is sufficiently large, then $f_{\lambda}(K)$ has exactly two unstable points, and if $\lambda$ is sufficiently small, then $f_{\lambda}(K)$ has exactly two stable points.

Similarly, a convex body $K$ with a significantly flattened or elongated \emph{bounding box} may not have equilibrium points but in some small regions of $\bd K$ including two small patches and a narrow ring as demonstrated by \cite{domokos2010pebbles}.
However, the number of equilibria within these areas can be arbitrary. Indeed, many objects have localized clusters of equilibria. For example, randomly generated fine polygonal/polyhedral discretizations of smooth surfaces often possess small groups of nearby equilibria, in place of isolated equilibria of the original surface.  Based on this observation, global equilibria consisting of clusters of local equilibria have been defined in \cite{domokos2010pebbles}. This concept was made precise in \cite{DLSz2012}, where the authors gave simple formulas for the numbers of stable, saddle and unstable points on a fine polyhedral approximation of a $C^2$-class convex surface $S$ in $\Re^3$ induced by an equidistant partition of the parameter range, based on the principal curvatures of the equilibrium points and their distances from the reference points. They called these quantities the \emph{imaginary equilibrium indices} of the equilibrium points of $S$, and carried out a similar investigation of $C^2$-class curves in $\Re^2$. In \cite{DLS2022}, the authors removed some technical assumptions from the results in \cite{DLSz2012}, and proved a similar result for the number of equilibrium points of a polygonal curve generated by randomly chosen points on the curve.

\subsection{The number of equilibria with respect to a general reference point}\label{subsec:inhom}

Let $K$ be a convex body in $\Re^d$. Let the function $n_K(c):\Re^d \to \mathbb{N}$ be defined as the number of equilibrium points with respect to the point $c$. We note that if $K$ has $C^2$-class boundary and strictly positive Gaussian curvature, then this quantity is equal to the number of points of $\bd(K)$ through which the line normal to $\bd (K)$ contains $c$. Thus, in the literature, questions related to this function are often referred to as \emph{problems of concurrent normals}. In this subsection we discuss results related to $n_K(\cdot)$.

\begin{figure}
    \centering
    \includegraphics[height=5 cm]{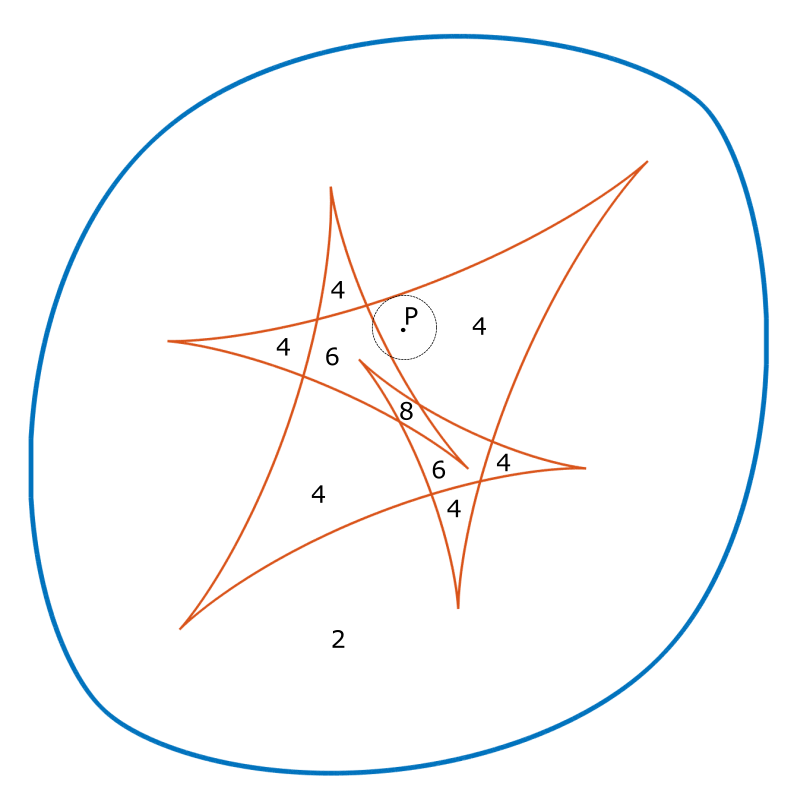}
    \includegraphics[height=5 cm]{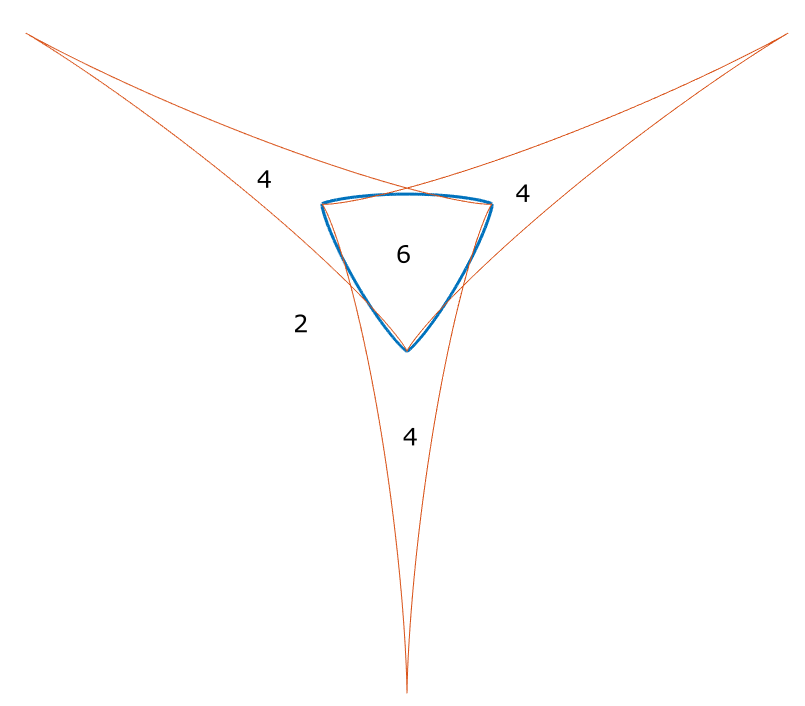}
    \caption{Two examples of smooth, planar curves with their evolutes, and the corresponding number of equilibria. In the left panel, the dashed circle illustrates the internal robustness of a point P. The right one is an example with $n_{min,K}=4$.}
    \label{fig:evolute}
\end{figure}

First, we consider plane convex bodies.
The function $n_K(\cdot)$ was first investigated by Apollonius of Perga \cite{apollonius1896apollonius} in the special case of an ellipse, and found that its value in a generic point is either 2 or 4.
It is known that for any $K$ in $\Re^2$ with $C^2$-class boundary and strictly positive curvature, $n_K(\cdot)$ is a piecewise constant function. Even more, it was proved by Allemann, Hungerb{\"u}hler and Wasem \cite{allemann2021equilibria} that $n_K(c)=2+2|m|$  where $m$ is the winding number of the evolute of $\bd(K)$ with respect to $c$. This result makes it possible to have a nice illustration of the function $n_K(\cdot)$: we decompose the plane into connected regions by the evolute of $\bd(K)$ (which in the following we call the evolute of $K$) and observe that within each region $n_K(\cdot)$ is a constant (Figure \ref{fig:evolute}). As a consequence, $n_K(c)$ is discontinuous at the points of the evolute, and $K$ is degenerate only with respect to the points of its evolute. 

It is trivial that the number of equilibrium points of any convex body $K$ with respect to any point $c$ is at least two; this follows from the observation that the points of $\bd (K)$ closest to and farthest from $c$ are equilibrium points. On the other hand, for any convex body with $C^2$-class boundary and strictly positive Gaussian curvature, $n_K(c)=2$ for all points $c$ sufficiently far from $K$. However, due to motivation by mechanics, it is often assumed that the reference point $c$ lies in the interior of the convex body $K$. Thus, we define the quantity
\[
n_{min,K}=\min \{ n_K(c) : c \in \inter (K)\}.
\]
We observe that there are plane convex bodies $K$ satisfying $n_{min,K} \geq 4$. A simple example for this is a smoothened regular triangle (Figure \ref{fig:evolute}).
Even higher values of $n_{min,K}$ are possible if we drop the differentiability properties for $\bd(K)$: for example, rectangles have eight equilibrium points with respect to all interior points or more generally, bricks in $\Re^d$ have $3^d-1$ equilibrium points with respect to all their interior points. The highest possible value of $n_{min,K}$ among polytopes and among smooth bodies appear to be unsolved problems. For example, it is an intriguing question if there exist a plane convex body $K$ with $C^2$-class boundary and strictly positive curvature that satisfies $n_{min,K}>4$. A related result by Grebennikov et al. \cite{grebennikov2022note} is that almost every normal through a boundary point to such a convex body in any dimensions contains a point $c$ with $n_K(c)\geq 6$.

In contrast, the quantity
$$
n_{max,K}=\max \{ n_K(c) : c \in \inter (K)\}
$$ 
has been studied by several authors. More specifically, it has been proven under mild smoothness assumptions that $n_{max,K}\geq 2d$ for any convex body with $C^1$-class boundary in $\Re^d$ with $2 \leq d \leq 4$ \cite{heil1985concurrent,pardon2012concurrent}, and the same statement in arbitrary dimensions has been known for long as \emph{concurrent normal conjecture}. 
In addition, the inequality $n_{max,K}\geq 6$ has been proven for any convex body $K$ in $\Re^d$ with an arbitrary value of $d$ \cite{heil1985concurrent}. In a recent paper another lower bound for $n_{max,K}$ over the family of $d$-dimensional convex polytopes was given by Nasonov, Panina and Siersma \cite{nasonov2024concurrent}. Namely, they proved the inequality $n_{max,K}\geq 2d+2$, and conjectured that the inequality $n_{max,K} \geq 10$ holds for any $3$-dimensional convex polytope.

Several works have addressed the average number of equilibria
\begin{align}   
\overline{n}_K=\frac{1}{\vol_d(K)} \int_{\inter(K)}n_K(x)dx
\end{align}
in a convex body $K$, for certain families of convex bodies.
A list of bounds, established in the papers 
\cite{domokos2015average,hann1993average,hug1995mean,santalo1944note,hann1996s,dumitrascu1998every}, have been collected in Table \ref{tab:average}.


\begin{table}[htb]
    \centering
    \begin{tabular}{|r||c| c | c | c |}
    \hline
       Dimension  & 2 & 3 & $d$  \\
       \hline\hline 
         $\mathcal{P}$ & $\stackrel{?}{>}4$ \cite{dumitrascu1998every}, $\leq 12$ \cite{hann1993average}  & $\stackrel{*}{>} 2$ \cite{domokos2015average} &$\stackrel{*}{>}2$\cite{domokos2015average} \\
        \hline 
        $\mathcal{F}$ &  $\stackrel{*}{\geq} 2$, & $\stackrel{*}{\geq} 2$  & $\stackrel{*}{\geq} 2$ \\
        &  $\leq 12$ \cite{hann1993average,hug1995mean}; $\stackrel{?}{\leq} 8$ 
        \cite{hann1996s} & $\leq 62$; $\stackrel{?}{\leq} 26$ 
        \cite{hann1996s}  & $\leq \frac{3}{2}^d \binom{2d}{d}-1$ \\
				\hline 
        $\mathcal{F}_{cs}$  & $\stackrel{*}{\leq} 8$ \cite{hann1993average,hug1995mean} & $\stackrel{*}{\leq} 26$ \cite{hann1993average,hug1995mean} & $\stackrel{*}{\leq} 3^d-1$ \cite{hann1993average,hug1995mean} \\
				\hline 
           $\mathcal{F}_{cc}$ & $\leq 6$ \cite{hann1993average,hug1995mean}  & $\leq 20$ \cite{domokos2015average} & $\leq \binom{2d}{d}$ \cite{domokos2015average} \\
					\hline 
        $\mathcal{F}_{cw}$ & $\stackrel{*}{\leq} \frac{2\pi}{\pi-\sqrt{3}}$ \cite{santalo1944note, chakerian1984number} & $\leq 20$ \cite{domokos2015average} & $\leq \binom{2d}{d}$ \cite{domokos2015average} \\
        \hline
    \end{tabular}
    \caption{A summary of known bounds of the average number of equilibria. Star means bounds are verified to be sharp, and question marks denote published conjectures. The notation in the first column is as follows: $\mathcal{P}$ and $\mathcal{F}$ mean the sets of convex polytopes, and convex bodies with $C^2$-class boundary and strictly positive curvature, respectively, in the appropriate dimension. For $x \in \{cs, cc, cw \}$, $\mathcal{F}_x$ denotes the subset of $\F$ consisting of centrally symmetric convex bodies, convex bodies containing all their centers of curvature, and bodies of constant width bodies. The bounds established by \cite{hug1995mean} do not require $C^2$-class boundary.
    }
    \label{tab:average}
\end{table}

\subsection{Equilibrium points of convex polyhedra}\label{subsec:polyhedra}

In modern times, the study of equilibrium points of convex polyhedra was probably started with a problem of Conway and Guy \cite{conway1966stability} in 1966 who conjectured that there is no monostable tetrahedron but there is a monostable convex polyhedron in $\Re^3$. These two questions were answered by Goldberg and Guy in \cite{goldberg1969stability} in 1969 (Figure \ref{fig:GombocGuy}). For a more detailed proof of the first problem, see also \cite{dawson1985monostatic}. In addition, in \cite{goldberg1969stability} Guy presented some problems regarding monostable polyhedra, stating that three of them are due to Conway (for similar statements in the literature, see e.g. \cite{CFG, dawson1999shape, DT2020}). These three questions appear also in the problem collection of Croft, Falconer and Guy \cite{CFG} as Problem B12. These problems were:
\begin{enumerate}
\item Can a monostable polyhedron in the Euclidean $3$-space $\Re^3$ have an $n$-fold axis of symmetry for $n > 2$?
\item What is the smallest possible ratio of diameter to girth for a monostable polyhedron?
\item What is the set of convex bodies uniformly approximable by monostable polyhedra, and does this contain the sphere?
\end{enumerate}
Related to these problems we recall that the \emph{girth} of a convex body in $\Re^3$ is the minimum perimeter of an orthogonal projection of the body onto a plane \cite{DT2020}. It is worth noting that, as stated in \cite{goldberg1969stability}, Conway showed that no body of revolution can be monostable, and also that the polyhedron constructed in \cite{goldberg1969stability} has a $2$-fold rotational symmetry.

The first two questions have recently been answered by L\'angi \cite{Langi_2022}, but the third one is still open. The main tool of the proofs of the result in \cite{Langi_2022} is a theorem stating that a nondegenerate convex body, under some mild condition on its curvatures, can be approximated arbitrarily well by a nondegenerate convex polyhedron with the same numbers of equilibria and with the same symmetry group. This result can be regarded as the counterpoint of the one in \cite{DLSz2012} about approximations of a convex surface by equidistant partitions of its parameter range, mentioned in Subsection~\ref{subsec:location}.

A problem proposed in \cite{goldberg1969stability} asks about the minimum dimension, if it exists, in which a $d$-simplex can be monostable.
This problem has been investigated by Dawson et al. \cite{dawson1985monostatic, dawson1998monostatic, dawson2001monostatic, dawson1999shape}, who proved that there is no monostable $d$-simplex if $d \leq 8$ and there is a monostable $11$-simplex. With regard to Problem XVI in \cite{Shephard1968}, asking about the minimum number of faces of a monostable polyhedron in $\Re^3$, the original construction of Guy \cite{goldberg1969stability} with $19$ faces (attributed also to Conway) was modified by Bezdek to obtain a monostable polyhedron with $18$ faces \cite{bezdek2011on}, while a computer-aided search by Reshetov \cite{reshetov2014unistable} yields a monostable polyhedron with $14$ faces. The non-existence of mono-unstable homogeneous tetrahedra has been proved in \cite{domokos2018balancing} where a construction of a mono-unstable polyhedron was also given.

We have seen that all equilibrium classes $\{ S,U \}$ contain nondegenerate convex polyhedra \cite{Langi_2022} and all classes $\{ S \}$  contain convex polygons \cite{Langi_2022}, as one can construct a fine polyhedral appoximation of a smooth convex body from the same class. On the other hand, these representatives have much more faces, edges, and vertices than the number of equilibria. Recall that every vertex, edge or face of a nondegenerate convex polyhedron contains at most one equilibrium point in its relative interior. Thus, an interesting related problem is to find the `simplest' convex polyhedron in each class in terms of the number of vertices, edges or faces. This leads to the following concept.

For an arbitrary nondegenerate, convex polyhedron $P$ the total number of faces, edges, and vertices without equilibria is called the \emph{complexity} of $P$, denoted by $C(P)$. Furthermore, for any $S,U \geq 1$, we denote the smallest complexity among polyhedra in the class $\{S,U\}_p$ by $C(S,U)$, and call it the \emph{complexity of equilibrium class $\{S,U\}_p$ }. In addition, recall a classical theorem of Steinitz \cite{Steinitz1, Steinitz2} stating that a triple of integers $(v,e,f)$ are the numbers of the vertices, edges and faces of a convex polyhedron, respectively, if, and only if $f \geq 4$, $\frac{f}{2}+2 \leq v \leq 2f-4$, and $e=v+f-2$. Based on this result, let us call a pair of integers $(f,v)$ a \emph{polyhedral pair} if they satisfy these inequalities, and let $\mathcal{PP}$ denote the set of all polyhedral pairs.

It has been shown in \cite{domokos2018balancing} that for all equilibrium classes with $S,U>1$, 
\[
C(S,U)=2\min \{ f+v-S-U : (f,v)\in \mathcal{PP} \}.
\]
In particular, this result implies that $C(S,U)=0$ if and only if $(S,U)$ is a polyhedral pair, or in other words, any primary equilibrium class $\{ S,U \}_p$, where $(S,U)$ is a polyhedral pair, can be represented by a \emph{minimal polyhedron}, that is, by a polyhedron that contains an equilibrium point in the relative interior of every face, edge and vertex. An illustration of polyhedral pairs, and some examples of minimal polyhedra are given in Figure \ref{fig:table}. There is no similar result for secondary and tertiary equilibrium classes, though the problem seems to be related to the problem of centering Koebe polyhedra (see \cite{domokos2018balancing}, and for more information on centering Koebe polyhedra, see \cite{Langi2021}).

While the above result answers our question for all classes except for monostable and mono-unstable ones, exact values of the complexities of the classes of monostatic polyhedra are unknown. In particular, the monostable polyhedron $P_{CG}$ of Conway and Guy in \cite{conway1966stability} is in class $ \{1,4 \}_p$ and has complexity $C(P_{CG})=96$. The ones $P_B$ and $P_R$, constructed by Bezdek and Reshetov, are in classes $\{ 1,3 \}_p$ and $\{ 1,2 \}_p$, respectively, having complexity $C(P_B)=64$ and $C(P_R)=70$. The authors of \cite{domokos2018balancing} extended these constructions for all classes $\{ 1, U \}_p$ with $U \geq 2$, and constructed polyhedra in classes $\{S, 1 \}_p$ for all $S \geq 2$, yielding the estimates $C(2,1) \leq 66$, $C(3,1) \leq 64$, and for all $S,U \geq 4$, $C(S,1) \leq 61 + 2 \left\lfloor \frac{S}{2} \right\rfloor$ and $C(1,U) \leq 92 + 2 \left\lfloor \frac{U}{2} \right\rfloor$. A lower bound $v\geq 7$ for the number $v$ of vertices of mono-unstable polyhedra was established by \cite{bozoki2024smallest} using a computational approach based on semi-definite programming. Presently there is no explicit estimate for the complexity of the class $ \{1,1 \}_p$. However it has been conjectured that $C(1,1)$ is at least approximately $1000$ \cite{varkonyi2006static}.

The picture outlined above changes significantly, if the centroid of the polyhedron is replaced by the center of mass of the vertices. In this case,
mono-unstable polyhedra with less than eight vertices do not exist \cite{bozoki2022mono}, but a mono-monostatic polyhedron with $21$ faces and $21$ vertices has been found \cite{domokos2023conway}. If the assumption of homogeneity is relaxed, and the reference point is permitted to be an arbitrary interior point of the polyhedron, then even some (non-regular) tetrahedra may be monostable as first pointed out by Conway (see \cite{dawson1999shape}) whereas other tetrahedra may be mono-unstable \cite{almadi2023equilibria}. It is also shown in \cite{almadi2023equilibria} that all monostable and mono-unstable tetrahedra belong to the equilibrium class $\{ 1,2 \}_p$ and $\{2,1\}_p$, respectively. In higher dimensions, even some regular polytopes may be monostable \cite{dawson1999shape}.

\section{Related problems and applications}\label{sec:appl}

The concepts of equilibrium and stability are closely related to the equilibria of floating bodies, which is the main focus of a review paper of Alfonseca et al. in preparation. The equilibrium positions of a convex body $K$ floating in a liquid correspond to the equilibrium positions in the sense of Definition \ref{defn:equilibrium} of a related body $U_V(K)$, called the \emph{Ulam floating body} of $K$ \cite{huang2019ulam}. The boundary of this body is also referred to as the \emph{surface of buoyancy} (Alfonseca et al., in preparation). For any $0< V < \vol_d(K)$, the surface of buoyancy contains the centroids of the intersections of $K$ with all closed half spaces $H$ that satisfy the condition $\vol_d(K \cap H) = V$. Thus, if $K$ is symmetric and $V = \frac{1}{2} \vol_d(K)$, then the Ulam floating body of $K$  coincides with the centroid body of $K$ (see Subsection~\ref{subsec:BP}). Moreover, in the limit cases $V \to 0$ and $V \to \vol_d(K)$, up to a homothety, $U_V(K)$ tends to the original body $K$. 

A possible generalization of the notion of an equilibrium point of a convex body is the following \cite{allemann2021equilibria}: given an angle $0 < \alpha < \frac{\pi}{2}$, we say that a point $q$ in the boundary of a convex body $K$ is an \emph{equilibrium point of $K$ of parameter $\alpha$}, with respect to $c$, if $K$ has a supporting hyperplane at $q$ with angle equal to $\frac{\pi}{2}-\alpha$ to $[c,q]$.
This definition is motivated by physical objects that can be balanced at a point on planes with slope angle $-\frac{\pi}{2} < \alpha < \frac{\pi}{2}$, provided that there is sufficient friction to prevent downhill slip. The extension of stability is also straightforward in two dimensions in the case $c=o$ by noting that such points correspond to critical points of the function $u \mapsto \cos\alpha h_K(u)+\sin\alpha s(u)$, $u \in \Sph^1$, where $s(u)$ is the (signed) arclength of the arc of $\bd(K)$ between $\rho_K(u) u \in \bd(K)$ and an arbitrary fixed reference point in $\bd(K)$. This function corresponds to the potential energy function of a heavy plane body rolling on a slope without slip.  Allemann, Hungerb{\"u}hler and Wasem \cite{allemann2021equilibria} showed that the number of equilibria of parameter $\alpha$ of a plane convex body $K$ is a piecewise constant function of $\alpha$, which drops to zero outside a bounded interval (with 0 in its interior) of the slope angle $\alpha$.
In dimensions $d \geq 3$ there is no straightforward way of defining stability, as strictly convex bodies can always escape from an equilibrium by means of a spin motion.

The concepts outlined in this paper can also be extended to equilibria of objects enclosed in a sphere of radius $\alpha$  \cite{varkonyi2016sensorless}; in this context the equilibrium points in Definition~\ref{defn:equilibrium} correspond to the limit case $\alpha\to \infty$. The equilibria of a convex body $K$ in a sphere of radius $\alpha$ are in one-to-one correspondence with the equilibrium points of the $\alpha$-hull of $K$, defined as the intersection of all balls of radius $\alpha$ that contain $K$ (for the definition of $\alpha$-convex sets and a systematic investigation of these bodies, the reader is referred to \cite{Mayer1935} and \cite{BLNP2007}, respectively). The number of equilibrium points of a plane convex body $K$ with respect to an interior point $c$, defined in this way, is a piecewise constant, strictly increasing function of $\alpha$, with minimum value of two for every $c$ different from the center of the smallest ball containing $K$. In case of $3$-dimensional convex bodies, this function is not necessarily monotonic, but it also drops to two if $\alpha$ is sufficiently close to the circumradius of the body. In other words, almost all objects become mono-monostatic inside a sufficiently small sphere \cite{varkonyi2016sensorless}.  

A further generalization considers equilibria of objects with multiple contact points on an underlying surface with arbitrary geometry. A partial characterization of the positions of the centers of mass corresponding to equilibria has been given by \cite{or2006computation,rimon2008general}, and rich geometric structures have been uncovered. The stability analysis of these systems relies heavily on the theory of dynamical systems \cite{posa2015stability,varkonyi2017lyapunov} and is beyond the scope of this work. 

We finish the paper by briefly mentioning a few  applications of the above results outside mathematics. The theorem of Winternitz (see Subsection~\ref{subsec:Grunbaum}) has recently found an application to data depth in statistics \cite{NSW2019}. Equilibrium classes (see Subsections \ref{subsec:primary}, \ref{subsec:sectert}) are used as a shape classification scheme in geology \cite{korvin2024shape,domokos2010pebbles,ludmany2023pebbles}. Flocks of equilibria (see Subsection~\ref{subsec:location}) help to understand fascinating geological formations called \emph{balancing rocks} \cite{bell1998dating,ludmany2023new}. 

Monostable objects play special roles in industrial object manipulation processes called part feeding  \cite{berkowitz1996designing}, and the statistical behaviour of randomly dropped objects with multiple stable equilibria also attracts much interest \cite{goldberg1999part,varkonyi2014estimating}. Monostable shapes also appear in nature \cite{domokos2008geometry,chiari2017self} as well as in robotics \cite{tunstel1999evolution,zhao2013msu}, and they allow spontaneous self-righting if an object topples. A further medical application of monostable objects concerns autonomous drug delivery  \cite{abramson2019ingestible}.

\bibliographystyle{plain}
\bibliography{literature}

\end{document}